# CONSTRUCTION OF NESTED SPACE-FILLING DESIGNS

BY PETER Z. G. QIAN[1], MINGYAO AI[2] AND C. F. JEFF WU[3]

*University of Wisconsin–Madison, Peking University and Georgia Institute of Technology*

New types of designs called *nested space-filling designs* have been proposed for conducting multiple computer experiments with different levels of accuracy. In this article, we develop several approaches to constructing such designs. The development of these methods also leads to the introduction of several new discrete mathematics concepts, including nested orthogonal arrays and nested difference matrices.

**1. Introduction.** Computer models are widely used in business, engineering and sciences to study complex real-world systems. The corresponding physical experimentation might otherwise be time-consuming, costly or even infeasible to conduct. Space-filling designs [Fang, Li and Sudjianto (2006) and Santner, Williams and Notz (2003)] have been widely used for conducting computer experiments. They include Latin hypercube designs [McKay, Conover and Beckman (1979)] and their improvements and variants [Butler (2001), Owen (1992, 1994b), Steinberg and Lin (2006), Tang (1993, 1998) and Ye (1998)]. Statistical properties of such designs have been studied in Loh (1996a, 1996b, 2008), Owen (1994a) and Stein (1987). Other types of space-filling designs are uniform designs [Fang et al. (2000)], quasi-Monte Carlo sequences [Niederreiter (1992)] and designs with uniform coverage [Dalal and Mallows (1998) and Lam, Welch and Young (2002)].

A large computer code, like a finite element analysis model, is often run at variable degrees of sophistication, resulting in multiple computer exper-

Received June 2008; revised November 2008.
[1]Supported by NSF Grant DMS-07-05206 and a faculty award from IBM.
[2]Supported by NNSF of China Grant 10671007 and NBRP of China Grant 2007CB512605.
[3]Supported by NSF Grant DMS-07-05261.

iments with different levels of accuracy and varying computational times. In this article, we consider the situation in which two such experiments are available, and one source is generally more accurate than the other but also more expensive to run. As in Qian and Wu (2008), the two experiments considered are called the *high-accuracy experiment* (HE) and *low-accuracy experiment* (LE). The problem of modeling data from HE and LE has attracted a recent surge of interests. Related work includes Goldstein and Rougier (2004), Higdon et al. (2004), Kennedy and O'Hagan (2000, 2001), Reese et al. (2004), Qian et al. (2006) and Qian and Wu (2008), among others. Most of these methods are based on flexible Gaussian process models [Fang, Li and Sudjianto (2006), Sacks et al. (1989), Santner, Williams and Notz (2003) and Welch et al. (1992)].

The sets of design points for LE and HE are denoted by $D_l$ and $D_h$. Throughout the paper, LE and HE are assumed to share the same set of factors and the design region, for both $D_l$ and $D_h$, are assumed to be the unit hypercube. As a suitable choice for $D_l$ and $D_h$, the notion of *nested space-filling designs* (NSFDs) was introduced in Qian, Tang and Wu (2009) (referred to as QTW hereinafter). The basic idea is to construct a special orthogonal array $A_1$ and use it to obtain an OA-based Latin hypercube design [Tang (1993)] for $D_l$. Take $A_2$ to be a subset of $A_1$ that becomes an orthogonal array itself after some level-collapsing, then obtain $D_h$ as the subarray of $D_l$ corresponding to $A_2$. The constructed $D_l$ and $D_h$ achieve low-dimensional uniformity. The nested relationship $D_h \subset D_l$ is appealing, which is also adopted in Kennedy and O'Hagan (2000), Qian et al. (2006) and Qian and Wu (2008). It implies that the size of $D_h$ is smaller than that of $D_l$ which is desirable because LE is cheaper than HE, and more LE runs can be afforded. From the modeling standpoint, this structure ensures that for every point in $D_h$, the outputs from both HE and LE are available, thus making it easier to model the differences of outputs between the two sources, and perform model adjustment.

We call the above special orthogonal array *nested orthogonal array* (NOA). Its formal definition will be given in the next section. A family of NOAs with fixed levels was constructed in QTW based on the Rao–Hamming method which will be reviewed in Section 2.3. In this article, we propose a new approach to constructing such arrays. The principal idea is to first construct nested difference matrices and then take the Kronecker product of a nested difference matrix and a standard orthogonal array to obtain an NOA. This method is motivated by the fact that constructing a nested difference matrix is probably easier than the direct construction of its corresponding NOA. Similar considerations have been used in constructing orthogonal arrays from difference matrices [Hedayat, Sloane and Stufken (1999), referred to as HSS hereinafter]. As a modification of this approach, we provide another method that uses existing NOAs to obtain new ones. These methods can produce many new NOAs and therefore new NSFDs. Several approaches for constructing NOAs with mixed levels will also be discussed.

The remainder of the article will unfold as follows. In Section 2, some notation and definitions are introduced. In Section 3, an approach based on multiplication tables of Galois fields to constructing nested difference matrices is proposed. In Section 4, a general approach to constructing NOAs with Kronecker product is presented. In Section 5, a method is introduced for constructing new NOAs from existing ones. In Section 6, construction of NOAs with nonprime power number of levels is considered. Construction of NOAs with mixed levels is given in Section 7. In Section 8, the problem of using NOAs to obtain NSFDs is discussed. Some discussions and concluding remarks are provided in Section 9.

## 2. Notation and definitions.

2.1. *Preliminaries.* Let $A = (a_{ij})$ be a Latin hypercube of $n$ runs for $m$ factors that is an $n \times m$ matrix where each column is a permutation of $1, \ldots, n$. Following McKay, Conover and Beckman (1979), a Latin hypercube design of $n$ runs in $m$ factors in the unit cube $[0,1)^m$ is generated through $x_{ij} = (a_{ij} - u_{ij})/n, 1 \le i \le n, 1 \le j \le m$, where $u_{ij}$'s are independent $U(0,1]$ random variables, and the $n$ design points are given by $(x_{i1}, \ldots, x_{im})$, $i = 1, \ldots, n$. When such a design is projected onto each of the $m$ factors, one and only one of the $n$ points falls within each of the $n$ small intervals defined by $[0, 1/n), [1/n, 2/n), \ldots, [(n-1)/n, 1)$.

A symmetrical orthogonal array (OA) of size $n$, $m$ constraints, $s$ levels, and strength $t \ge 2$, is an $n \times m$ matrix with entries from a set of $s$ levels, usually taken as $1, \ldots, s$, such that for every $n \times t$ submatrix, the $s^t$ level combinations occurs equally often. Regular fractional factorial designs, as discussed in Wu and Hamada (2000), are the most familiar examples of orthogonal arrays. In the article, we consider only OAs with strength two, denoted by $OA(n, m, s)$. Asymmetrical OAs will be discussed in Section 7.

Let $A$ be an $OA(n, m, s)$ with its $s$ levels denoted by $1, \ldots, s$. Then in every column of $A$, each level occurs $q = n/s$ times. For each column of $A$, if we replace the $q$ ones by a permutation of $1, \ldots, q$, replace the $q$ twos by a permutation of $q+1, \ldots, 2q$, and so on, we obtain an OA-based Latin hypercube [Tang (1993)]. In addition to achieving maximum stratification in one dimension, OA-based Latin hypercubes have attractive space-filling properties when projected onto 2 dimensions.

A difference matrix (DM) is a $b \times c$ array with entries from a finite abelian group $\mathcal{A}$ with $g$ elements, such that every element of $\mathcal{A}$ appears equally often in the vector difference between any two columns of the array [Bose and Bush (1952)]. We will denote such an array by $D(b, c, g)$. If $\mathcal{A}$ is the additive group associated with a Galois field, we simply say its elements come from the associated field. For any $D(b,c,g)$, a column is defined to be *uniform* in $\mathcal{A}$ if it contains each element of $\mathcal{A}$ equally often. By subtracting the first column from all columns, any $D(b,c,g)$ can always be converted to a difference matrix of the form

$$[0_b \quad D^{(0)}], \tag{1}$$

where $0_b$ is the $b$-dimensional zero vector and every column of $D^{(0)}$ is uniform in $\mathcal{A}$.

Let $A = (a_{ij})$ and $B = (b_{ij})$ be, respectively, $m \times n$ and $u \times v$ matrices with entries from an abelian group $\mathcal{A}$ with binary operation $*$ (usually addition or multiplication). The *Kronecker product* of $A$ and $B$ [Shrikhande (1964)], denoted by $A \otimes B$, is defined to be the $mu \times nv$ matrix

$$A \otimes B = \begin{bmatrix} a_{11} * B & \cdots & a_{1n} * B \\ \vdots & & \vdots \\ a_{m1} * B & \cdots & a_{mn} * B \end{bmatrix},$$

where $a_{ij} * B$ denotes the $u \times v$ matrix with entries $a_{ij} * b_{rs}$, $1 \leq r \leq u, 1 \leq s \leq v$. Throughout this article $*$ always denotes addition.

2.2. *Galois field projections.* For every prime $p$ and every integer $u \geq 1$, there exists a Galois field (or finite field) $\mathrm{GF}(p^u)$ of order $p^u$. The additive group $\mathrm{GF}(p^u)$ is cyclic, and the multiplicative group $\mathrm{GF}(p^u)/\{0\}$ is cyclic, allowing easy calculations under multiplication. Throughout, the elements of any Galois field or any subset of a Galois field are arranged in lexicographical order.

Unless stated otherwise, let $s_1 = p^{u_1}$ and $s_2 = p^{u_2}$ be powers of the same prime $p$ with integers $u_1 > u_2 \geq 1$. Throughout, let $F$ denote $\mathrm{GF}(s_1)$ with an irreducible polynomial $p_1(x)$, and $G$ denote $\mathrm{GF}(s_2)$ with an irreducible polynomial $p_2(x)$. Let $f(x)$ denote the elements of $F$ and $g(x)$ the elements of $G$, respectively. In condensed notation, let $\alpha_0, \ldots, \alpha_{s_1-1}$ denote the elements of $F$ and $\beta_0, \ldots, \beta_{s_2-1}$ the elements of $G$ with $\alpha_0 = 0$ and $\beta_0 = 0$. Next, we discuss two projections from $F$ to $G$, serving as a basis for later development.

The first projection, denoted by $\phi$, is taken from Bose and Bush (1952). For any $f(x) = a_0 + a_1 x + \cdots + a_{u_2-1} x^{u_2-1} + \cdots + a_{u_1-1} x^{u_1-1} \in F$, $\phi(f(x))$ is defined by

$$\phi(f(x)) = a_0 + a_1 x + \cdots + a_{u_2-1} x^{u_2-1}. \tag{2}$$

Because $\phi$ works by truncating all $x$ powers of degree $u_2$ or higher, we call it the *truncation projection*.

The second projection, denoted by $\varphi$, is proposed in QTW. For any $f(x) \in F$, $\varphi(f(x))$ is defined by

$$\varphi(f(x)) = f(x) (\mathrm{mod}\, p_2(x)). \tag{3}$$

Because $\varphi$ works by taking modulus residues, we call it the *modulus projection*.

EXAMPLE 1. Let $p = 2$, $u_1 = 3$ and $u_2 = 2$, giving $s_1 = 8$ and $s_2 = 4$. Use $p_1(x) = x^3 + x + 1$ for GF(8) and $p_2(x) = x^2 + x + 1$ for GF(4). Then the projection $\varphi$ is given as $\{0, x^2 + x + 1\} \to 0$, $\{1, x^2 + x\} \to 1$, $\{x, x^2 + 1\} \to x$, $\{x + 1, x^2\} \to x + 1$.

Let $\delta$ be either of the projections described above. For an array $D$ with entries from $F$, $\delta(D)$ denotes the array obtained from $D$ after the levels of its entries are collapsed according to $\delta$. Clearly, the entries of $\delta(D)$ take values in $G$.

Notice that for any $\alpha_i, \alpha_j \in F$,

$$\delta(\alpha_i + \alpha_j) = \delta(\alpha_i) + \delta(\alpha_j). \tag{4}$$

This means that the two operations $\delta$ and $+$ are interchangeable, which is critical to the constructions in Sections 4, 5 and 7.

2.3. *Nested space-filling designs and nested orthogonal arrays.* Now we give a formal definition of NOAs, which underly the construction of NSFDs in QTW. Let $A_1$ be an OA$(n_1, k, s_1)$. Suppose there is a subarray of $A_1$ with size $n_2$, denoted by $A_2$, and there is a projection $\delta$ that collapses the $s_1$ levels of $A_1$ into $s_2$ levels. Further suppose $A_2$ becomes an OA$(n_2, k, s_2)$ after the levels of its entries are collapsed according to $\delta$. Then $A_1$, or more precisely $(A_1, A_2)$, is an NOA, denoted by NOA$(A_1, A_2)$ or NOA$(A_1, A_2, \delta)$. To be emphatic about a small OA being nested within a larger OA, we say $A_1$ "contains" $\delta(A_2)$.

Let $(A_1, A_2)$ be an NOA defined above. Construction of an NSFD is done as follows. The array $A_1$ is used to generate an OA-based Latin hypercube design $D_l$. Let $D_h$ denote the subset of $D_l$ corresponding to $A_2$. Then $D_l$, or more precisely $(D_l, D_h)$, is an NSFD, where both $D_l$ and $D_h$ achieve uniformity in low dimensions.

The family of NOA$(A_1, A_2)$, constructed in QTW by using the Rao–Hamming method, has the following set of parameters:

(i) $A_1$ is an OA$(n_1, m_2, s_1)$, where $n_1 = s_1^k$, $m_2 = (s_2^k - 1)/(s_2 - 1)$ and $k \geq 2$ is an integer;

(ii) $A_2$ is a subarray of $A_1$ and $\varphi(A_2)$ is an OA$(n_2, m_2, s_2)$ with $n_2 = s_2^k$.

This construction works for $2u_2 \leq u_1 + 1$.

2.4. *Nested difference matrices.* Let $D_1$ be a $D(b_1, c, s_1)$ with entries from $F$. Suppose there is a subarray of $D_1$ with $b_2$ rows denoted by $D_2$, and a projection $\delta$ that collapses the $s_1$ levels of $D_1$ into the $s_2$ levels of $G$. Further suppose $D_2$ is a $D(b_2, c, s_2)$ if the levels of its entries are collapsed according to $\delta$. Then $D_1$, or more precisely $(D_1, D_2)$, is called a *nested difference matrix* (NDM), denoted by $\mathrm{NDM}(D_1, D_2)$ or $\mathrm{NDM}(D_1, D_2, \delta)$. To be emphatic about a smaller DM being nested within a larger DM, we say $D_1$ "contains" $\delta(D_2)$.

**3. Construction of nested difference matrices.** In this section, we propose an approach based on multiplication tables of Galois fields to constructing NDMs. It works for any $u_1 > u_2 \geq 1$. Here the projection $\phi$ in (2) is used. For a scalar $a$ and a vector $c = (c_1, \ldots, c_m)'$, $a+c$ denotes $(a+c_1, \ldots, a+c_m)'$, where $'$ stands for vector transpose. Similarly, $a + A$ denotes the elementwise sum of a scalar $a$ and a matrix $A$. We focus on the case of $p = 2$ and briefly discuss the case of $p = 3$ in the end of the section. Two sets or vectors are defined to be *disjoint* if they have no element in common. Because the constructions in Section 4 can use a small NDM and a standard OA to generate a larger NOA, here we construct NDMs with up to 16 columns. Throughout, we use the irreducible polynomial $p(x) = x^u + x + 1$ for any $\mathrm{GF}(2^u)$, $u \geq 1$. Unless stated otherwise, let $r_{-1} = (0)$, $r_0 = (0, 1)'$, $r_m = (0, 1, x, x+1, \ldots, x^m + \cdots + x + 1)'$, $m \geq 1$. Note that $r_m$ has $2^{m+1}$ elements.

A $D(s_1, s_1, s_1)$ can be obtained by constructing the $s_1 \times s_1$ multiplication table of $\mathrm{GF}(s_1)$, where the rows and columns are labeled by all distinct elements of $\mathrm{GF}(s_1)$. Hereinafter, in describing such a table, we call a row (or column) labeled with an element $f(x) \in \mathrm{GF}(s_1)$ as "row (or column) $f(x)$."

3.1. *A $D(2^{m+1}, 2^2, 2^{m+1})$ containing a $D(2^m, 2^2, 2^m)$ with $m \geq 2$.* Let $F = \mathrm{GF}(2^{u_1})$ and $G = \mathrm{GF}(2^{u_2})$ with $u_1 = m+1$, $u_2 = m$ and $m \geq 2$. Let $D_0$ be the multiplication table of $F$. By taking columns $r_1$ of $D_0$, obtain a matrix $D_1$.

Collect the elements of $F$ into two vectors:

(5) $$g_1 = (r'_{m-2}, x^{m-1} + r'_{m-2})' \quad \text{and} \quad g_2 = x^m + g_1,$$

where the $i$th element in $g_2$ equals its counterpart in $g_1$ plus $x^m$. Now place the rows of $D_1$ in two clusters: the top one comprising those labeled with $r_{m-2}$ and $x^m + r_{m-2}$, and the bottom one with $x^{m-1} + r_{m-2}$ and $x^m + x^{m-1} + r_{m-2}$. This arrangement may look abstract at this moment but will become clear after Theorem 1. Table 1 gives $\phi(D_1)$, where, for $m = 2$, the entries need to be taken modulo $p_1(x) = x^{u_1} + x + 1$ and then collapsed according to $\phi$.

Take $D_2$ to be the submatrix of $D_1$ consisting of the rows labeled with $r_{m-2}$ and $x^m + x^{m-1} + r_{m-2}$. Because $r_{m-2}$ is the set of polynomials of order at most $m-2$, $r_{m-2}$ and $x^{m-1} + r_{m-2}$ are disjoint and their union is $GF(2^m)$. The following is a simple result regarding columns $x$ and $x+1$ of $\phi(D_1)$.

LEMMA 1. (i) *The vectors* $(x+1)r_{m-2}$ *and* $x^{m-1} + (x+1)r_{m-2}$ *are disjoint and their union is* $\mathrm{GF}(2^m)$;

(ii) *the vectors* $(x+1)r_{m-2}$ *and* $(x^{m-1} + x + 1) + (x+1)r_{m-2}$ *are disjoint and their union is* $\mathrm{GF}(2^m)$.

(iii) *the vectors* $xr_{m-2}$ *and* $(x+1) + xr_{m-2}$ *are disjoint and their union is* $GF(2^m)$.

PROOF. (i) It suffices to show that $(x+1)r_{m-2}$ and $x^{m-1} + (x+1)r_{m-2}$ are disjoint. Assuming the contrary, then there are two elements $\alpha_1$ and $\alpha_2$ from $r_{m-2}$ such that $(x+1)\alpha_1 = x^{m-1} + (x+1)\alpha_2$, implying $(x+1)(\alpha_1 - \alpha_2) = x^{m-1}$. This is impossible because $x+1$ does not divide $x^{m-1}$.

(ii) It follows from (i) by noting that $(x+1) + (x+1)r_{m-2}$ has the same set of elements as $(x+1)r_{m-2}$.

(iii) Assuming the contrary, then there are two elements $\alpha_1$ and $\alpha_2$ from $r_{m-2}$ such that $\alpha_1 - \alpha_2 - 1 = x^{-1}$, a contradiction. □

THEOREM 1. *Consider* $D_1$ *and* $D_2$ *constructed above. For* $m \geq 2$, *we have:*

(i) *the matrix* $D_1$ *is a* $D(2^{m+1}, 2^2, 2^{m+1})$;
(ii) *the matrix* $\phi(D_2)$ *is a* $D(2^m, 2^2, 2^m)$.

PROOF. Only (ii) needs a proof. Because the elements $\{0, 1, x, x+1\}$, used to label the columns of $D_1$, form an additive group, it suffices to show that columns $1, x, x+1$ of $\phi(D_2)$ are uniform in $GF(2^m)$. Note that, due to the grouping scheme in (5), columns 1 and $x$ of $\phi(D_2)$ in Table 1 are exactly an half fraction of those of $\phi(D_1)$. Then it remains to show that column $x+1$ of $\phi(D_2)$ is uniform in $\mathrm{GF}(2^m)$. This follows from Lemma 1 as $(x+1)+(x+1)r_{m-2}$ and $(x+1)r_{m-2}$ have the same set of elements. □

TABLE 1
*The matrix $\phi(D_1)$ obtained from $D_1$ in Theorem 1*

|  | **0** | **1** | **$x$** | **$x+1$** |
|---|---|---|---|---|
| $r_{m-2}$ | 0 | $r_{m-2}$ | $xr_{m-2}$ | $(x+1)r_{m-2}$ |
| $x^m + r_{m-2}$ | 0 | $r_{m-2}$ | $(x+1) + xr_{m-2}$ | $(x+1) + (x+1)r_{m-2}$ |
| $x^{m-1} + r_{m-2}$ | 0 | $x^{m-1} + r_{m-2}$ | $xr_{m-2}$ | $x^{m-1} + (x+1)r_{m-2}$ |
| $x^m + x^{m-1} + r_{m-2}$ | 0 | $x^{m-1} + r_{m-2}$ | $(x+1) + xr_{m-2}$ | $(x^{m-1} + x + 1) + (x+1)r_{m-2}$ |

EXAMPLE 2 [A $D(2^2, 2, 2^2)$ containing a $D(2, 2, 2)$]. Although this example has only two columns, we include it here because its construction is similar to those in Theorem 1. Let $F = \mathrm{GF}(2^2)$ and $G = \mathrm{GF}(2)$. Let $D_0$ be the multiplication table of $F$ given by

|       | **0** | **1** | **$x$** | **$x+1$** |
|-------|-------|-------|---------|-----------|
| 0     | 0     | 0     | 0       | 0         |
| 1     | 0     | 1     | $x$     | $x+1$     |
| $x$   | 0     | $x$   | $x+1$   | 1         |
| $x+1$ | 0     | $x+1$ | 1       | $x$       |

Take $D_1$ be the first two columns of $D_0$. Obtain $D_2$ as the submatrix of $D_1$ consisting of rows 0 and 1. The matrix $D_1$ is a $D(2^2, 2, 2^2)$, and $\phi(D_2)$ is a $D(2, 2, 2)$ given by

$$\begin{bmatrix} 0 & 0 \\ 0 & 1 \end{bmatrix}.$$

EXAMPLE 3 [A $D(2^3, 2^2, 2^3)$ containing a $D(2^2, 2^2, 2^2)$]. Let $F = \mathrm{GF}(2^3)$ and $G = \mathrm{GF}(2^2)$. Take $D_1$ to be the columns of the multiplication table of $F$ labeled with $r_1$ given by

|             | **0** | **1**       | **$x$**     | **$x+1$**   |
|-------------|-------|-------------|-------------|-------------|
| 0           | 0     | 0           | 0           | 0           |
| 1           | 0     | 1           | $x$         | $x+1$       |
| $x^2$       | 0     | $x^2$       | $x+1$       | $x^2+x+1$   |
| $x^2+1$     | 0     | $x^2+1$     | 1           | $x^2$       |
| $x$         | 0     | $x$         | $x^2$       | $x^2+x$     |
| $x+1$       | 0     | $x+1$       | $x^2+x$     | $x^2+1$     |
| $x^2+x$     | 0     | $x^2+x$     | $x^2+x+1$   | 1           |
| $x^2+x+1$   | 0     | $x^2+x+1$   | $x^2+1$     | $x$         |

From Theorem 1, $D_1$ is a $D(2^3, 2^2, 2^3)$, $D_2$ is the submatrix of $D_1$ consisting of rows $0, 1, x^2+x, x^2+x+1$, and $\phi(D_2)$ is a $D(2^2, 2^2, 2^2)$ given by

|             | **0** | **1**   | **$x$**   | **$x+1$** |
|-------------|-------|---------|-----------|-----------|
| 0           | 0     | 0       | 0         | 0         |
| 1           | 0     | 1       | $x$       | $x+1$     |
| $x^2+x$     | 0     | $x$     | $x+1$     | 1         |
| $x^2+x+1$   | 0     | $x+1$   | 1         | $x$       |

3.2. *A $D(2^{m+2}, 2^2, 2^{m+2})$ containing a $D(2^m, 2^2, 2^m)$ with $m \geq 2$.* Let $F = \mathrm{GF}(2^{u_1})$ and $G = \mathrm{GF}(2^{u_2})$ with $u_1 = m+2$, $u_2 = m$ and $m \geq 2$. Let $D_0$ be the multiplication table of $F$. By taking columns $r_1$ of $D_0$, obtain a matrix $D_1$.

Collect the elements of $F$ into two vectors:

(6)
$$g_1 = (r'_{m-2}, x^{m-1} + r'_{m-2}, x^{m+1} + r'_{m-2}, x^{m+1} + x^{m-1} + r'_{m-2})' \quad \text{and}$$
$$g_2 = x^m + g_1.$$

Now place the rows of $D_1$ in four clusters. From top to bottom, their row labels are: cluster 1 with $r_{m-2}$ and $x^m + r_{m-2}$; cluster 2 with $x^{m-1} + r_{m-2}$ and $x^m + x^{m-1} + r_{m-2}$; cluster 3 with $x^{m+1} + r_{m-2}$ and $x^{m+1} + x^m + r_{m-2}$; and cluster 4 with $x^{m+1} + x^{m-1} + r_{m-2}$ and $x^{m+1} + x^m + x^{m-1} + r_{m-2}$. Table 2 gives $\phi(D_1)$, where for $m = 2$ or $3$, the entries need to be taken modulo $p_1(x) = x^{u_1} + x + 1$ and then collapsed according to $\phi$.

Take $D_2$ to be the submatrix of $D_1$ consisting of the rows labeled with $r_{m-2}$ and $x^{m+1} + x^m + x^{m-1} + r_{m-2}$.

THEOREM 2. *Consider $D_1$ and $D_2$ constructed above. For $m \geq 2$, we have:*

(i) *the matrix $D_1$ is a $D(2^{m+2}, 2^2, 2^{m+2})$;*
(ii) *the matrix $\phi(D_2)$ is a $D(2^m, 2^2, 2^m)$.*

The proof of this theorem is similar to that of Theorem 1 and therefore omitted.

EXAMPLE 4 [A $D(2^4, 2^2, 2^4)$ containing a $D(2^2, 2^2, 2^2)$]. Let $F = \mathrm{GF}(2^4)$ and $G = \mathrm{GF}(2^2)$. From Theorem 2, $D_1$ is a $D(2^4, 2^2, 2^4)$, $D_2$ is the submatrix of $D_1$ consisting of the rows labeled with $(0, 1, x^3 + x^2 + x, x^3 + x^2 + x + 1)'$,

TABLE 2
*The matrix $\phi(D_1)$ obtained from $D_1$ in Theorem 2*

| | **0** | **1** | **$x$** | **$x+1$** |
|---|---|---|---|---|
| $r_{m-2}$ | 0 | $r_{m-2}$ | $xr_{m-2}$ | $(x+1)r_{m-2}$ |
| $x^m + r_{m-2}$ | 0 | $r_{m-2}$ | $xr_{m-2}$ | $(x+1)r_{m-2}$ |
| $x^{m-1} + r_{m-2}$ | $0\,x^{m-1} + r_{m-2}$ | $xr_{m-2}$ | $x^{m-1} + (x+1)r_{m-2}$ | |
| $x^m + x^{m-1} + r_{m-2}$ | $0\,x^{m-1} + r_{m-2}$ | $xr_{m-2}$ | $x^{m-1} + (x+1)r_{m-2}$ | |
| $x^{m+1} + r_{m-2}$ | 0 | $r_{m-2}$ | $(x+1) + xr_{m-2}$ | $(x+1) + (x+1)r_{m-2}$ |
| $x^{m+1} + x^m + r_{m-2}$ | 0 | $r_{m-2}$ | $(x+1) + xr_{m-2}$ | $(x+1) + (x+1)r_{m-2}$ |
| $x^{m+1} + x^{m-1} + r_{m-2}$ | $0\,x^{m-1} + r_{m-2}$ | $(x+1) + xr_{m-2}$ | $(x^{m-1} + x + 1) + (x+1)r_{m-2}$ | |
| $x^{m+1} + x^m + x^{m-1} + r_{m-2}$ | $0\,x^{m-1} + r_{m-2}$ | $(x+1) + xr_{m-2}$ | $(x^{m-1} + x + 1) + (x+1)r_{m-2}$ | |

and $\phi(D_2)$ is a $D(2^2, 2^2, 2^2)$ given by

|  | 0 | 1 | $x$ | $x+1$ |
|---|---|---|---|---|
| 0 | 0 | 0 | 0 | 0 |
| 1 | 0 | 1 | $x$ | $x+1$ |
| $x^3 + x^2 + x$ | 0 | $x$ | $x+1$ | 1 |
| $x^3 + x^2 + x + 1$ | 0 | $x+1$ | 1 | $x$ |

3.3. *A $D(2^{m+2}, 2^3, 2^{m+2})$ containing a $D(2^{m+1}, 2^3, 2^m)$ with $m \geq 2$.* Let $F = \mathrm{GF}(2^{u_1})$ and $G = \mathrm{GF}(2^{u_2})$ with $u_1 = m+2$, $u_2 = m$ and $m \geq 2$. Let $D_0$ denote the multiplication table of $F$. By taking columns $r_2$ of $D_0$, obtain a matrix $D_1$.

Collect the elements of $F$ into two vectors:

(7)
$$g_1 = (r'_{m-2}, x^{m-1} + r'_{m-2}, x^{m+1} + r'_{m-2}, x^{m+1} + x^{m-1} + r'_{m-2})' \quad \text{and}$$
$$g_2 = x^m + g_1.$$

Now place the rows of $D_1$ in four clusters. From top to bottom, their row labels are: cluster 1 with $r_{m-2}$ and $x^m + r_{m-2}$; cluster 2 with $x^{m-1} + r_{m-2}$ and $x^m + x^{m-1} + r_{m-2}$; cluster 3 with $x^{m+1} + r_{m-2}$ and $x^{m+1} + x^m + r_{m-2}$; and cluster 4 with $x^{m+1} + x^{m-1} + r_{m-2}$ and $x^{m+1} + x^m + x^{m-1} + r_{m-2}$. Table 3 gives columns $x^2 + r_1$ of $\phi(D_1)$, where, for $m = 2$ or 3, the entries need to be taken modulus $p_1(x)$ and then collapsed according to $\phi$, and

$$\alpha_1 = x^2(r'_{m-3}, r'_{m-3})',$$
$$\alpha_2 = ((x^2 + 1)r'_{m-3}, x^{m-2} + (x^2 + 1)r'_{m-3})',$$
$$\alpha_3 = ((x^2 + x)r'_{m-3}, x^{m-1} + (x^2 + x)r'_{m-3})',$$
$$\alpha_4 = ((x^2 + x + 1)r'_{m-3}, x^{m-1} + x^{m-2} + (x^2 + x + 1)r'_{m-3})'.$$

Take $D_2$ to be the submatrix of $D_1$ consisting of rows $r_{m-2}$, $x^m + x^{m-1} + r_{m-2}$, $x^{m+1} + r_{m-2}$, $x^{m+1} + x^m + x^{m-1} + r_{m-2}$.

THEOREM 3. *Consider $D_1$ and $D_2$ constructed above. For $m \geq 2$, we have:*

(i) *the matrix $D_1$ is a $D(2^{m+2}, 2^3, 2^{m+2})$;*
(ii) *the matrix $\phi(D_2)$ is a $D(2^{m+1}, 2^3, 2^m)$.*

Its proof is similar to that of Theorem 1 and therefore omitted.

EXAMPLE 5 [A $D(2^5, 2^3, 2^5)$ containing a $D(2^4, 2^3, 2^3)$]. Let $F = \mathrm{GF}(2^5)$ and $G = \mathrm{GF}(2^3)$. We have $g_1 = (0, 1, x, x+1, x^2, x^2+1, x^2+x, x^2+x+1, x^4, x^4+1, x^4+x, x^4+x+1, x^4+x^2, x^4+x^2+1, x^4+x^2+x, x^4+x^2+x+1)'$ and $\alpha_1 = (0, x^2, 0, x^2)'$, $\alpha_2 = (0, x^2+1, x, x^2+x+1)'$, $\alpha_3 = (0, x^2+x, x^2, x)'$,

TABLE 3
*Columns $x^2 + r_1$ of $\phi(D_1)$ obtained from $D_1$ in Theorem 3*

| | $x^2$ | $x^2+1$ | $x^2+x$ | $x^2+x+1$ |
|---|---|---|---|---|
| $r_{m-2}$ | $\alpha_1$ | $\alpha_2$ | $\alpha_3$ | $\alpha_4$ |
| $x^m + r_{m-2}$ | $(x+1)+\alpha_1$ | $(x+1)+\alpha_2$ | $(x+1)+\alpha_3$ | $(x+1)+\alpha_4$ |
| $x^{m-1} + r_{m-2}$ | $\alpha_1$ | $x^{m-1}+\alpha_2$ | $\alpha_3$ | $x^{m-1}+\alpha_4$ |
| $x^m + x^{m-1} + r_{m-2}$ | $(x+1)+\alpha_1$ | $(x^{m-1}+x+1)+\alpha_2$ | $(x+1)+\alpha_3$ | $(x^{m-1}+x+1)+\alpha_4$ |
| $x^{m+1} + r_{m-2}$ | $(x^2+x)+\alpha_1$ | $(x^2+x)+\alpha_2$ | $(x^2+1)+\alpha_3$ | $(x^2+1)+\alpha_4$ |
| $x^{m+1} + x^m + r_{m-2}$ | $(x^2+1)+\alpha_1$ | $(x^2+1)+\alpha_2$ | $(x^2+x)+\alpha_3$ | $(x^2+x)+\alpha_4$ |
| $x^{m+1} + x^{m-1} + r_{m-2}$ | $(x^2+x)+\alpha_1$ | $(x^{m-1}+x^2+x)+\alpha_2$ | $(x^2+1)+\alpha_3$ | $(x^{m-1}+x^2+1)+\alpha_4$ |
| $x^{m+1} + x^m + x^{m-1} + r_{m-2}$ | $(x^2+1)+\alpha_1$ | $(x^{m-1}+x^2+x)+\alpha_2$ | $(x^2+x)+\alpha_3$ | $(x^{m-1}+x^2+1)+\alpha_4$ |

$\alpha_4 = (0, x^2+x+1, x^2+x, 1)'$. From Theorem 3, $D_1$ is a $D(2^5, 2^3, 2^5)$, $D_2$ is the submatrix of $D_1$ consisting of rows $0, 1, x, x+1, x^3+x^2, x^3+x^2+1, x^3+x^2+x, x^3+x^2+x+1, x^4, x^4+1, x^4+x, x^4+x+1, x^4+x^3+x^2, x^4+x^3+x^2+1, x^4+x^3+x^2+x, x^4+x^3+x^2+x+1$, and $\phi(D_2)$ is a $D(2^4, 2^3, 2^3)$, with columns $x^2 + r_1$ of $\phi(D_2)$ given by

|                       | $x^2$         | $x^2+1$       | $x^2+x$       | $x^2+x+1$     |
|-----------------------|---------------|---------------|---------------|---------------|
| $0$                   | $0$           | $0$           | $0$           | $0$           |
| $1$                   | $x^2$         | $x^2+1$       | $x^2+x$       | $x^2+x+1$     |
| $x$                   | $0$           | $x$           | $x^2$         | $x^2+x$       |
| $x+1$                 | $x^2$         | $x^2+x+1$     | $x$           | $1$           |
| $x^3+x^2$             | $x+1$         | $x^2+x+1$     | $x+1$         | $x^2+x+1$     |
| $x^3+x^2+1$           | $x^2+x+1$     | $x$           | $x^2+1$       | $0$           |
| $x^3+x^2+x$           | $x+1$         | $x^2+1$       | $x^2+x+1$     | $1$           |
| $x^3+x^2+x+1$         | $x^2+x+1$     | $0$           | $1$           | $x^2+x$       |
| $x^4$                 | $x^2+x$       | $x^2+x$       | $x^2+1$       | $x^2+1$       |
| $x^4+1$               | $x$           | $x+1$         | $x+1$         | $x$           |
| $x^4+x$               | $x^2+x$       | $x^2$         | $1$           | $x+1$         |
| $x^4+x+1$             | $x$           | $1$           | $x^2+x+1$     | $x^2$         |
| $x^4+x^3+x^2$         | $x^2+1$       | $1$           | $x^2+x$       | $x$           |
| $x^4+x^3+x^2+1$       | $1$           | $x^2$         | $0$           | $x^2+1$       |
| $x^4+x^3++x^2+x$      | $x^2+1$       | $x+1$         | $x$           | $x^2$         |
| $x^4+x^3+x^2+x+1$     | $1$           | $x^2+x$       | $x^2$         | $x+1$         |

3.4. *Some extensions.* Some extensions of the proposed method are considered here. Similar to Sections 3.1–3.3, we can construct the following two families of NDMs: (a) a $D(2^{m+3}, 2^3, 2^{m+3})$ containing a $D(2^{m+1}, 2^3, 2^m)$ and (b) a $D(2^{m+3}, 2^4, 2^{m+3})$ containing a $D(2^{m+2}, 2^4, 2^m)$ with $m \geq 2$. For brevity we present the case with $m = 2$, where $F = \mathrm{GF}(2^5)$ and $G = \mathrm{GF}(2^2)$. By taking columns $r_3$ of the multiplication table of $F$, we obtain a matrix $D_1$. Clearly, $D_1$ is a $D(2^5, 2^4, 2^5)$.

Collect the elements of $F$ into

$$g_1 = (r_0', x + r_0', x^3 + r_0', x^3 + x + r_0', x^4 + r_0', x^4 + x + r_0',$$
$$x^4 + x^3 + r_0', x^4 + x^3 + x + r_0')'$$

and $g_2 = x^2 + g_1$. Note that, for the columns labeled with $r_2$, the $i$th row in $g_1$ is the same as its counterpart in $g_2$. Let $D_2$ be the submatrix of $D_1$ consisting of rows $(0, 1, x^3+x, x^3+x+1, x^4, x^4+1, x^4+x^3+x, x^4+x^3+x+1)$. It is easy to see that $\phi(D_2)$ is a $D(2^3, 2^3, 2^2)$ given by

|  | **0** | **1** | $\boldsymbol{x}$ | $\boldsymbol{x+1}$ | $\boldsymbol{x^2}$ | $\boldsymbol{x^2+1}$ | $\boldsymbol{x^2+x}$ | $\boldsymbol{x^2+x+1}$ |
|---|---|---|---|---|---|---|---|---|
| 0 | 0 | 0 | 0 | 0 | 0 | 0 | 0 | 0 |
| 1 | 0 | 1 | $x$ | $x+1$ | 0 | 1 | $x$ | $x+1$ |
| $x^3+x$ | 0 | $x$ | 0 | $x$ | $x+1$ | 1 | $x+1$ | 1 |
| $x^3+x+1$ | 0 | $x+1$ | $x$ | 1 | $x+1$ | 0 | 1 | $x$ |
| $x^4$ | 0 | 0 | $x+1$ | $x+1$ | $x$ | $x$ | 1 | 1 |
| $x^4+1$ | 0 | 1 | 1 | 0 | $x$ | $x+1$ | $x+1$ | $x$ |
| $x^4+x^3+x$ | 0 | $x$ | $x+1$ | 1 | 1 | $x+1$ | $x$ | 0 |
| $x^4+x^3+x+1$ | 0 | $x+1$ | 1 | $x$ | 1 | $x$ | 0 | $x+1$ |

Take $D_3$ to be the submatrix of $D_1$ consisting of rows $(0,1), x^2+(x, x+1), x^2+(x^3, x^3+1), (x^3+x, x^3+x+1), x^2+(x^4, x^4+1), (x^4+x, x^4+x+1), (x^4+x^3, x^4+x^3+1), x^2+(x^4+x^3+x, x^4+x^3+x+1)$ of $D_1$. It is easy to verify that $\phi(D_3)$ is a $D(2^4, 2^4, 2^2)$ with columns $x^3 + r_2$ of $\phi(D_3)$ given by

|  | $\boldsymbol{x^3}$ | $\boldsymbol{x^3+1}$ | $\boldsymbol{x^3+x}$ | $\boldsymbol{x^3+x+1}$ | $\boldsymbol{x^3+x^2}$ | $\boldsymbol{x^3+x^2+1}$ | $\boldsymbol{x^3+x^2+x}$ | $\boldsymbol{x^3+x^2+x+1}$ |
|---|---|---|---|---|---|---|---|---|
| 0 | 0 | 0 | 0 | 0 | 0 | 0 | 0 | 0 |
| 1 | 0 | 1 | $x$ | $x+1$ | 0 | 1 | $x$ | $x+1$ |
| $x^2+x$ | $x+1$ | 1 | $x+1$ | 1 | $x+1$ | 1 | $x+1$ | 1 |
| $x^2+x+1$ | $x+1$ | 0 | 1 | $x$ | $x+1$ | 0 | 1 | $x$ |
| $x^2+x^3$ | 1 | 1 | 1 | 1 | $x$ | $x$ | $x$ | $x$ |
| $x^2+x^3+1$ | 1 | 0 | $x+1$ | $x$ | $x$ | $x+1$ | 0 | 1 |
| $x^3+x$ | $x$ | 0 | $x$ | 0 | 1 | $x+1$ | 1 | $x+1$ |
| $x^3+x+1$ | $x$ | 1 | 0 | $x+1$ | 1 | $x$ | $x+1$ | 0 |
| $x^2+x^4$ | $x+1$ | $x+1$ | 0 | 0 | 1 | 1 | $x$ | $x$ |
| $x^2+x^4+1$ | $x+1$ | $x$ | $x$ | $x+1$ | 1 | 0 | 0 | 1 |
| $x^4+x$ | 0 | $x$ | $x+1$ | 1 | $x$ | 0 | 1 | $x+1$ |
| $x^4+x+1$ | 0 | $x+1$ | 1 | $x$ | $x$ | 1 | $x+1$ | 0 |
| $x^4+x^3$ | $x$ | $x$ | 1 | 1 | $x+1$ | $x+1$ | 0 | 0 |
| $x^4+x^3+1$ | $x$ | $x+1$ | $x+1$ | $x$ | $x+1$ | $x$ | $x$ | $x+1$ |
| $x^2+x^4$ $+x^3+x$ | 1 | $x+1$ | $x$ | 0 | 0 | $x$ | $x+1$ | 1 |
| $x^2+x^4+x^3$ $+x+1$ | 1 | $x$ | 0 | $x+1$ | 0 | $x+1$ | 1 | $x$ |

The proposed method can be extended to construct NDMs with $p = 3$. Note that the presentation of this extension is more involved because the irreducible polynomials for $GF(3^u)$, $u \geq 1$, do not have a unified form. [In contrast, we can use $p(x) = x^u + x + 1$ for any $GF(2^u)$, $u \geq 1$.] For brevity we provide examples from a useful family: a $D(3^{m+1}, 3^2, 3^{m+1})$ containing a $D(3^m, 3^2, 3^m)$ with $m \geq 1$.

Here let $r_{-1} = (0)$, $r_0 = (0, 1, 2)'$, $r_m = (0, 1, 2, x, x+1, x+2, \ldots, 2x^m, 2x^m + 1, \ldots, 2x^m + 2x^{m-1} + \cdots + 2x + 2)'$ with $m \geq 1$. Note that here $r_m$ has $3^{m+1}$ elements. Let $F = \text{GF}(3^{u_1})$ with an irreducible polynomial $p_1(x)$ and $G = \text{GF}(3^{u_2})$ with an irreducible polynomial $p_2(x)$, where $u_1 = m + 1$, $u_2 = m$ and $m \geq 1$. Let $D_0$ be the multiplication table of $F$. By taking columns $r_1$ of $D_0$, obtain a matrix $D_1$. Clearly, $D_1$ is a $D(3^{m+1}, 3^2, 3^{m+1})$.

Collect the elements of $F$ into three vectors:

$$
\begin{aligned}
g_1 &= (r'_{m-2}, x^{m-1} + r'_{m-2}, 2x^{m-1} + r'_{m-2})', \\
g_2 &= x^m + g_1 \quad \text{and} \quad g_3 = 2x^m + g_1.
\end{aligned}
\tag{8}
$$

As a consequence of this grouping scheme, for any column labeled with $r_0$ in $\phi(D_1)$, the rows labeled with $g_1$ are the same as those labeled with $g_2$ or those labeled with $g_3$. This convenient structure implies that we only need to focus on columns $x + r_0$ and $2x + r_0$ in the construction. The key is to find a subset of $D_1$ in which these columns are uniform in $G$. Some examples are given.

EXAMPLE 6 [A $D(3^3, 3^2, 3^3)$ containing a $D(3^2, 3^2, 3^2)$]. Let $F = \text{GF}(3^3)$ with $p_1(x) = x^3 + 2x + 1$ and $G = \text{GF}(3^2)$ with $p_2(x) = x^2 + x + 2$. Let $D_0$ be the multiplication table of $F$. By taking columns $r_1$ of $D_0$, obtain a matrix $D_1$ which is a $D(3^3, 3^2, 3^3)$. Let $D_2$ be the submatrix of $D_1$ consisting of rows $r_0, 2x^2 + x + r_0$ and $x^2 + 2x + r_0$. It is easy to see that $\phi(D_2)$ is a $D(2^3, 2^2, 2^3)$ with columns $x, x+1, x+2, 2x, 2x+1, 2x+2$ given by

|             | $x$    | $x+1$  | $x+2$  | $2x$   | $2x+1$ | $2x+2$ |
|-------------|--------|--------|--------|--------|--------|--------|
| 0           | 0      | 0      | 0      | 0      | 0      | 0      |
| 1           | $x$    | $x+1$  | $x+2$  | $2x$   | $2x+1$ | $2x+2$ |
| 2           | $2x$   | $2x+2$ | $2x+1$ | $x$    | $x+2$  | $x+1$  |
| $2x^2 + x$  | $2x+1$ | 1      | $x+1$  | $x+2$  | $2x+2$ | 2      |
| $2x^2+x+1$  | 1      | $x+2$  | $2x$   | 2      | $x$    | $2x+1$ |
| $2x^2+x+2$  | $x+1$  | $2x$   | 2      | $2x+2$ | 1      | $x$    |
| $x^2+2x$    | $x+2$  | 2      | $2x+2$ | $2x+1$ | $x+1$  | 1      |
| $x^2+2x+1$  | $2x+2$ | $x$    | 1      | $x+1$  | 2      | $2x$   |
| $x^2+2x+2$  | 2      | $2x+1$ | $x$    | 1      | $2x$   | $x+2$  |

EXAMPLE 7 [A $D(3^4, 3^2, 3^4)$ containing a $D(3^3, 3^2, 3^3)$]. Let $F = \text{GF}(3^4)$ with $p_1(x) = x^4 + x + 2$ and $G = \text{GF}(3^3)$ with $p_2(x) = x^3 + 2x + 1$. Let $D_0$ be the multiplication table of $F$. By taking columns $r_1$ of $D_0$, obtain a matrix $D_1$, which is a $D(3^4, 3^2, 3^4)$. Let $D_2$ be the submatrix of $D_1$ consisting of rows $r_1, 2x^3 + x^2 + r_1$ and $x^3 + 2x^2 + r_1$. Columns $x, x+1, x+2, 2x, 2x+$

$1, 2x+2$ of $\phi(D_2)$ are given by

|  | $x$ | $x+1$ | $x+2$ |
|---|---|---|---|
| $r_1$ | $xr_1$ | $(x+1)r_1$ | $(x+2)r_1$ |
| $2x^3+x^2+r_1$ | $(x+2)+xr_1$ | $(x^2+x+2)+(x+1)r_1$ | $(2x^2+x+2)+(x+2)r_1$ |
| $x^3+2x^2+r_1$ | $(2x+1)+xr_1$ | $(2x^2+2x+1)+(x+1)r_1$ | $(x^2+2x+1)+(x+2)r_1$ |
|  | $2x$ | $2x+1$ | $2x+2$ |
| $r_1$ | $2xr_1$ | $(2x+1)r_1$ | $(2x+2)r_1$ |
| $2x^3+x^2+r_1$ | $(2x+1)+2xr_1$ | $(x^2+2x+1)+(2x+1)r_1$ | $(2x^2+2x+1)+(2x+2)r_1$ |
| $x^3+2x^2+r_1$ | $(x+2)+2xr_1$ | $(2x^2+x+2)+(2x+1)r_1$ | $(x^2+x+2)+(2x+2)r_1$ |

Similar to Lemma 1 it is easy to show that any two of $(x+1)r_1$, $(x+2)r_1$ and $(2x+1)r_1$ are disjoint and the union of the three is $r_2$. Hence the columns of $\phi(D_2)$ are uniform in $r_2$ and $\phi(D_2)$ is a $D(3^3, 3^2, 3^3)$.

**4. Constructing nested orthogonal arrays with Kronecker products.** In this section we present a general approach to constructing NOAs. It generates an NOA by taking the Kronecker product of an NDM and a standard OA. Let $\delta$ be either $\phi$ in (2) or $\varphi$ in (3) unless stated otherwise.

The following lemma [Bose and Bush (1952)] says that taking the Kronecker product of an OA and a DM gives a larger OA.

LEMMA 2. *If $D$ is a $D(b, c, s)$ and $A$ is an $\mathrm{OA}(n, k, s)$, and both are based on the same abelian group $\mathcal{A}$, then the array $H = A \otimes D$ is an $\mathrm{OA}(nb, kc, s)$.*

For $\delta(A \otimes D)$, we have:

LEMMA 3. *If $D$ is a $D(b, c, s)$ and $A$ is an $\mathrm{OA}(n, k, s)$, and both are based on $\mathrm{GF}(s)$, then*

$$(9) \qquad \delta(A \otimes D) = \delta(A) \otimes \delta(D).$$

This lemma can be readily verified by using the definition of $\delta$ and (4). It basically says the two operations $\delta$ and $\otimes$ in (9) are interchangeable, which is key to the constructions to be proposed later.

Now recall a classical result from Addleman and Kempthorne (1961).

LEMMA 4. *If a factor in an OA has $s_1$ levels and $s_2|s_1$, then it can be replaced by a new factor with $s_2$ levels by partitioning the $s_1$ symbols into $s_2$ groups of size $s_1/s_2$ and by replacing the symbols in the same group with a common symbol. The resulting array is still an OA.*

Note that if the $s_1$ and $s_2$ levels in this lemma come from $\mathrm{GF}(s_1)$ and $\mathrm{GF}(s_2)$, respectively, the condition $s_2|s_1$ clearly holds and the required level collapsing can be done through using $\delta$.

Here is a similar result for difference matrices.

LEMMA 5. *If $D$ is a $D(b, c, s_1)$ based on $\mathrm{GF}(s_1)$, then $\delta(D)$ is a $D(b, c, s_2)$.*

This lemma can be readily proved by following the definitions of difference matrices and $\delta$.

Now we are ready to present the details of the proposed construction. Let $A$ be an $\mathrm{OA}(n, m, s_1)$ based on $\mathrm{GF}(s_1)$. Let $(D_1, D_2, \delta)$ be an NDM constructed in Section 3, where $D_1$ is a $D(b_1, c, s_1)$ based on $\mathrm{GF}(s_1)$ and $D_2$ is a submatrix of $D_1$, and $\delta(D_2)$ is a $D(b_2, c, s_2)$ based on $\mathrm{GF}(s_2)$. Put

$$(10) \qquad H_1 = A \otimes D_1 \quad \text{and} \quad H_2 = A \otimes D_2.$$

THEOREM 4. *For $H_1$ and $H_2$ in (10), the array $(H_1, H_2, \delta)$ is an NOA, where $H_1$ is an $\mathrm{OA}(nr_1, mc, s_1)$, $H_2$ is a submatrix of $H_1$ and $\delta(H_2)$ is an $\mathrm{OA}(nr_2, mc, s_2)$.*

This theorem can be readily verified by following Lemmas 2–5.

EXAMPLE 8. Let $p = 2$, $u_1 = 3$, $u_2 = 2$, giving $s_1 = 8$ and $s_2 = 4$. Take an $\mathrm{NDM}(D_1, D_2, \phi)$ from Example 3, where $D_1$ is a $D(8, 4, 8)$ and $\phi(D_2)$ is a $D(4, 4, 4)$. The projection $\phi$ is as follows: $\{0, x^2\} \to 0, \{1, x^2 + 1\} \to 1, \{x, x^2 + x\} \to x, \{x + 1, x^2 + x + 1\} \to x + 1$. Let $A$ be a trivial orthogonal array $\mathrm{OA}(8, 1, 8)$, the column vector listing all elements of $\mathrm{GF}(8)$. From Theorem 4, $H_1$ is an $\mathrm{OA}(64, 4, 8)$, $H_2$ is a submatrix of $H_1$, and $\phi(H_2)$ is an $\mathrm{OA}(32, 4, 4)$.

Note that the construction (10) is not restricted to use NDMs from Section 3. Here is an example.

EXAMPLE 9. Let $D_0$ be the $D(12, 12, 4)$ [Seberry (1979)] given in the Appendix. Take $D_1$ to be the submatrix of $D_0$ consisting of columns 1, 3, 4 and 5. It can be verified that $D_1$ is a $D(12, 4, 4)$. Take $D_2$ to be the submatrix of $D_1$ consisting of rows 1,2,4,5. Let $\delta$ be a projection by deleting the first digits of the entries in $D_2$. Clearly, $\delta(D_2)$ is a $D(4, 4, 2)$. Let $A$ be the $\mathrm{OA}(64, 21, 4)$ constructed by using the Rao–Hamming method (HSS). Put $H_1 = A \otimes D_1$ and $H_2 = A \otimes D_2$. Then the array $(H_1, H_2, \delta)$ is an NOA, where $H_1$ is an $\mathrm{OA}(768, 84, 4)$, $H_2$ is a submatrix of $H_1$ and $\delta(H_2)$ is an $\mathrm{OA}(256, 84, 2)$.

**5. Obtaining new nested orthogonal arrays from existing ones.** As a modification of the method in the previous section, we discuss here a procedure for obtaining new NOAs from existing ones. Let $(A_1, A_2, \varphi)$ be an arbitrary NOA constructed in QTW, where $A_1$ is an $OA(n_1, m, s_1)$, $A_2$ is a submatrix of $A_1$ and $\varphi(A_2)$ is an $OA(n_2, m, s_2)$. Let $D$ be a $D(b, c, s_1)$ based on $GF(s_1)$. Put

(11) $$H_1 = A_1 \otimes D \quad \text{and} \quad H_2 = A_2 \otimes D.$$

THEOREM 5. *For $H_1$ and $H_2$ in (11), we have:*

(i) *the matrix $H_1$ is an $OA(n_1 b, mc, s_1)$;*
(ii) *the matrix $H_2$ is a submatrix of $H_1$ and $\varphi(H_2)$ is an $OA(n_2 b, mc, s_2)$.*

This theorem can be readily verified by following Lemmas 2, 3 and 5.

EXAMPLE 10. Let $p = 2$, $u_1 = 3$, $u_2 = 2$, giving $s_1 = 8$ and $s_2 = 4$. We use $p_1(x) = x^3 + x + 1$ for $GF(8)$ and $p_2(x) = x^2 + x + 1$ for $GF(4)$. The condition $2u_2 \leq u_1 + 1$ is satisfied and the projection $\varphi$ is as follows. $\{0, x^2 + x + 1\} \to 0, \{1, x^2 + x\} \to 1, \{x, x^2 + 1\} \to x, \{x+1, x^2\} \to x + 1$. Take an NOA$(A_1, A_2, \varphi)$ from Section 2.3, where $A_1$ is an $OA(64, 5, 8)$ and $A_2$ is the following submatrix of $A_1$

$$\begin{bmatrix} 0 & 0 & 0 & 0 & 0 \\ 1 & 0 & 1 & x & x+1 \\ x & 0 & x & x^2 & x^2+x \\ x+1 & 0 & x+1 & x^2+x & x^2+1 \\ 0 & 1 & 1 & 1 & 1 \\ 1 & 1 & 0 & x+1 & x \\ x & 1 & x+1 & x^2+1 & x^2+x+1 \\ x+1 & 1 & x & x^2+x+1 & x^2 \\ 0 & x & x & x & x \\ 1 & x & x+1 & 0 & 1 \\ x & x & 0 & x^2+x & x^2 \\ x+1 & x & 1 & x^2 & x^2+x+1 \\ 0 & x+1 & x+1 & x+1 & x+1 \\ 1 & x+1 & x & 1 & 0 \\ x & x+1 & 1 & x^2+x+1 & x^2+1 \\ x+1 & x+1 & 0 & x^2+1 & x^2+x \end{bmatrix}.$$

We have $\varphi(A_2)$ is an OA$(16,5,4)$. Let $D$ be the multiplication table of GF(8). Then $D$ is a $D(2^3,2^3,2^3)$ and $\varphi(D)$ is a $D(2^3,2^3,2^2)$ given by

|           | 0 | 1   | $x$   | $x+1$ | $x^2$ | $x^2+1$ | $x^2+x$ | $x^2+x+1$ |
|-----------|---|-----|-------|-------|-------|---------|---------|-----------|
| 0         | 0 | 0   | 0     | 0     | 0     | 0       | 0       | 0         |
| 1         | 0 | 1   | $x$   | $x+1$ | 0     | 1       | $x$     | $x+1$     |
| $x$       | 0 | $x$ | 0     | $x$   | $x+1$ | 1       | $x+1$   | 1         |
| $x+1$     | 0 | $x+1$ | $x$ | 1     | $x+1$ | 0       | 1       | $x$       |
| $x^2$     | 0 | 0   | $x+1$ | $x+1$ | $x$   | $x$     | 1       | 1         |
| $x^2+1$   | 0 | 1   | 1     | 0     | $x$   | $x+1$   | $x+1$   | $x$       |
| $x^2+x$   | 0 | $x$ | $x+1$ | 1     | 1     | $x+1$   | $x$     | 0         |
| $x^2+x+1$ | 0 | $x+1$ | 1   | $x$   | 1     | $x$     | 0       | $x+1$     |

From Theorem 5, $H_1$ is an OA$(512, 40, 8)$, $H_2$ is a submatrix of $H_1$, and $\varphi(H_2)$ is an OA$(128, 40, 8)$.

**6. Construction of nested orthogonal arrays with nonprime power number of levels.** In this section, we construct NOAs with nonprime power number of levels. This construction complements the methods in the previous two sections, where NOAs with prime power number of levels are constructed. First we introduce a simple projection, denoted by $\rho_a$, for any integer $a \geq 1$, to be

(12) $$\rho_a(u) = u(\mathrm{mod}\, a).$$

The following lemma gives some properties of $\rho_a$:

LEMMA 6. (i) *If $a, b \geq 1$ are integers with $b|a$, then $\rho_b(\rho_a(u)) = \rho_b(u)$;*
(ii) *for any integer $a \geq 1$, $\rho_a(u_1 + u_2) = \rho_a(\rho_a(u_1) + \rho_a(u_2))$.*

We now use this projection to construct a family of NOAs based on the *zero-sum array* (HSS). For an integer $s$, let $\mathbb{Z}$ denote the residue classes modulo $s$. Let $s_1, s_2 \geq 1$ be integers with $s_2|s_1$. Let $F$ denote $\mathbb{Z}_{s_1}$ and $G$ denote $\mathbb{Z}_{s_2}$. Obtain an $s_1^2 \times 3$ matrix $A_1$, where the first two columns have each of the $s_1^2$ possible 2-tuples from $F \times F$ as a row, and for row $(i,j)$ in the first two columns, its corresponding entry in the third column is taken as $-(i+j)(\mathrm{mod}\, s_1)$. Take $A_2$ to be the submatrix of $A_1$ consisting of rows $(i,j)$, $0 \leq i, j \leq s_2 - 1$, in the first two columns.

THEOREM 6. *For $A_1$ and $A_2$ constructed above, we have:*

(i) *the matrix $A_1$ is an OA$(s_1^2, 3, s_1)$;*
(ii) *the matrix $A_2$ is a submatrix of $A_1$ and $\rho_{s_2}(A_2)$ is an OA$(s_2^2, 3, s_2)$.*

This theorem can be readily verified by following Lemma 6.

As a straightforward extension of Theorem 5, we can take the Kronecker product of an NOA from Theorem 6 and a standard DM to obtain a new

NOA. By extending Theorem 4, we can take the Kronecker product of an NDM with nonprime power number of levels and a standard OA to obtain an NOA. Here is an example.

EXAMPLE 11. Obtain a matrix $D_1$ by suppressing the first digits of all entries of the $D(12,6,12)$ in the Appendix. It can be verified that $D_1$ is a $D(12,6,6)$. Let $D_2$ be the submatrix of $D_1$ consisting of rows $1,4,5,6,8$ and $12$. Clearly, $\rho_3(D_2)$ is a $D(6,6,3)$. Let $A$ be the $OA(36,3,6)$ obtained by taking the first three columns of Table 7C.8 in Wu and Hamada (2000). Put $H_1 = A \otimes D_1$ and $H_2 = A \otimes D_2$. The array $(H_1, H_2, \rho_3)$ is an NOA, where $H_1$ is an $OA(432,18,6)$, $H_2$ is a submatrix of $H_1$ and $\rho_3(H_2)$ is an $OA(216,18,3)$.

**7. Construction of nested orthogonal arrays with mixed levels.** In this section, we discuss the issue of constructing NOAs with mixed levels. The key here is to embed nested structures in the constructions of asymmetrical (mixed) OAs, like those in Wang and Wu (1991) (referred to as WW hereinafter) and Wang (1996). Such an embedding can be done in various ways as described in the remainder of the section. We use $\mathrm{OA}(n, s_1^{\gamma_1} \cdots s_k^{\gamma_k})$ to denote an asymmetrical OA.

7.1. *Using nested orthogonal arrays and Wang–Wu method.* This construction makes use of the Kronecker products in (11) and the Wang–Wu method in WW. For $1 \leq j \leq v$, let $s_{j1}$ and $s_{j2}$ be powers of the same prime $p_j$ with integers $u_{j1} > u_{j2} \geq 1$. The primes $p_j$'s are assumed to be all distinct. Suppose $A_1$ is an $\mathrm{OA}(n_1, s_{11}^{k_1} \cdots s_{v1}^{k_v})$ and can be partitioned as

$$A_1 = [A_{11} \quad \cdots \quad A_{v1}],$$

where each $A_{j1}$ comes from an $\mathrm{NOA}(A_{j1}, A_{j2}, \delta_j)$, $A_{j1}$ is an $\mathrm{OA}(n_1, k_j, s_{j1})$ based on $\mathrm{GF}(s_{j1})$, $A_{j2}$ is a submatrix of $A_{j1}$ and $\delta_j(A_{j2})$ is an $\mathrm{OA}(n_2, k_j, s_{j2})$ based on $\mathrm{GF}(s_{j2})$.

For $1 \leq j \leq v$, let $D(j)$ denote a $D(b, c_j, s_{j1})$ with entries from $\mathrm{GF}(s_{j1})$. Put

$$H_1 = [A_{11} \otimes D(1) \cdots A_{v1} \otimes D(v) B_1]$$

and

$$H_2 = [A_{12} \otimes D(1) \cdots A_{v2} \otimes D(v) B_2],$$

where $C = (0, \ldots, b-1)'$, $B_1 = (C', \ldots, C')'$ represents a $b$-level factor with $C$ appearing $n_1$ times and $B_2 = (C', \ldots, C')'$ represents a $b$-level factor with $C$ appearing $n_2$ times.

THEOREM 7. *For $H_1$ and $H_2$ constructed above, we have:*

(i) *the matrix $H_1$ is an* $\mathrm{OA}(bn_1, b^1 s_{11}^{k_1 c_1} \cdots s_{v1}^{k_v c_v})$;

(ii) *the matrix $H_2$ is a submatrix of $H_1$ and $H_2$ is an* $\mathrm{OA}(bn_2, b^1 s_{12}^{k_1 c_1} \cdots s_{v2}^{k_v c_v})$ *after the levels of the $s_{j1}$-level factors are collapsed according to $\delta_j$ for $j = 1, \ldots, v$.*

This theorem can be readily verified by following the result on the generalized Kronecker product in WW and the definition of $\delta_j$.

EXAMPLE 12 [An $\mathrm{OA}(288, 6^6 4^{12})$ containing an $\mathrm{OA}(72, 3^6 2^{12})$]. Let $A_1$ be an $\mathrm{OA}(24, 6^1 4^1)$ formed by taking all level combinations of a factor at six levels, $0, 1, 2, 3, 4, 5$, and a factor at four levels, $00, 01, 10, 11$. Take $A_2$ to be the subarray of $A_1$ consisting of all level combinations of $0, 1, 2$ and $00, 01$. Let $D_1$ be the $D(12, 6, 6)$ and $D_2$ the $D(12, 12, 4)$ from the Appendix. Put $H_1 = [A_{11} \otimes D_1, A_{21} \otimes D_2]$ and $H_2 = [A_{12} \otimes D_1, A_{22} \otimes D_2]$. From Theorem 7, $H_1$ is an $\mathrm{OA}(288, 6^6 4^{12})$ and $H_2$ becomes an $\mathrm{OA}(72, 3^6 2^{12})$ after the following level collapsing: for the 6-level factors, using $\{0, 3\} \to 0$, $\{1, 4\} \to 1$, $\{2, 5\} \to 2$; for the 4-level factors, deleting the first digit and retaining the second, for example, both 01 and 11 are projected to 1.

7.2. *Using nested difference matrices and Wang–Wu method.* This construction makes use of the Kronecker products in (10) and the Wang–Wu method in WW.

For $1 \leq j \leq v$, let $s_{j1}$ and $s_{j2}$ be powers of the same prime $p_j$ with integers $u_{j1} > u_{j2} \geq 1$. The primes $p_j$'s are assumed to be all different. Suppose $A$ is an $\mathrm{OA}(n, s_{11}^{k_1} \cdots s_{v1}^{k_v})$ and can be partitioned as

$$A = [A_1 \ \cdots \ A_v],$$

where $A_j$ is an $\mathrm{OA}(n, k_j, s_{j1})$ based on $\mathrm{GF}(s_{j1})$. Let $D$ be a partitioned matrix

$$[D_1(1) \ \cdots \ D_1(v)],$$

where $D_1(j)$ comes from an $\mathrm{NDM}(D_1(j), D_2(j), \delta_j)$, $D_1(j)$ is a $D(b_1, c_j, s_{j1})$ based on $\mathrm{GF}(s_{j1})$, $D_2(j)$ is a submatrix of $D_1(j)$, and $\delta_j(D_2(j))$ is a $D(b_2, c_j, s_{j2})$ based on $\mathrm{GF}(s_{j2})$.

Put

$$H_1 = [A_1 \otimes D_1(1) \cdots A_v \otimes D_1(v) B_1]$$

and

$$H_2 = [A_1 \otimes D_2(1) \cdots A_v \otimes D_2(v) B_2],$$

where $C_1 = (0, \ldots, b_1 - 1)'$, $B_1 = (C_1', \ldots, C_1')'$ represents a $b_1$-level factor with $C_1$ appearing $n$ times, $C_2 = (0, \ldots, b_2 - 1)'$ is a subvector of $C_1$ with $b_2$ elements and $B_2 = (C_2', \ldots, C_2')'$ represents a $b_2$-level factor with $C_2$ appearing $n$ times.

THEOREM 8. *For $H_1$ and $H_2$ constructed above, we have:*

(i) *the matrix $H_1$ is an $\mathrm{OA}(b_1 n, b_1^1 s_{11}^{k_1 c_1} \cdots s_{v1}^{k_v c_v})$;*

(ii) *the matrix $H_2$ is a submatrix of $H_1$ and $H_2$ becomes an $\mathrm{OA}(b_2 n, b_2^1 s_{12}^{k_1 c_1} \cdots s_{v2}^{k_v c_v})$ after the levels of the $s_{j1}$-level factors are collapsed according to $\delta_j$ for $j = 1, \ldots, v$.*

This theorem can be readily verified by following the result on the generalized Kronecker product in WW and the definition of $\delta_j$.

7.3. *Using a nested nonorthogonal mixed matrix and a special mixed difference matrix.* Wang (1996) constructs an asymmetrical OA using a mixed DM and a nonorthogonal matrix with mixed levels. Unlike the Wang–Wu method, this construction does not use OAs and therefore can give asymmetrical OAs with more flexible run sizes. Here we modify it to construct NOAs with mixed levels.

For $j = 1, 2$, let $s_{j1}$ and $s_{j2}$ be powers of the same prime $p_j$ with integers $u_{j1} > u_{j2} \geq 1$. For $j = 1, 2$, choose an $\mathrm{NDM}(D_{j1}, D_{j2}, \delta_j)$, where $D_{j1}$ is a $D(n_1, k_j, s_{j1})$ based on $\mathrm{GF}(s_{j1})$, $D_{j2}$ is a submatrix of $D_{j1}$ and $\delta_j(D_{j2})$ is a $D(n_2, k_j, s_{j2})$ with entries from $\mathrm{GF}(s_{j2})$. Construct an $\mathrm{NDM}(D_{01}, D_{02}, \delta_0)$, where $\delta_0 = \delta_1 \times \delta_2$, $D_{01}$ is a $D(n_1, k_0, s_{11} s_{21})$ based on $\mathrm{GF}(s_{11}) \times \mathrm{GF}(s_{21})$, $D_{02}$ is a submatrix of $D_{01}$ and $\delta_0(D_{02})$ is a $D(n_2, k_0, s_{12} s_{22})$ based on $\mathrm{GF}(s_{12}) \times \mathrm{GF}(s_{22})$. For $j = 1, 2$, let $\sigma_j(\cdot)$ denote the operation of taking the $j$th component of every entry in a matrix whose entries are represented by two digits. For $j = 1, 2$, further assume the augmented matrix $[\sigma_j(D_{01}), D_{j1}]$ is a $D(n_1, k_0 + k_j, s_{j1})$ and $\delta_j[\sigma_j(D_{02}), D_{j2}]$ is a $D(n_2, k_0 + k_j, s_{j2})$. For $j = 1, 2$, let $C_j$ be the column vector comprising all level combinations of $\mathrm{GF}(s_{1j})$ and $\mathrm{GF}(s_{2j})$.

Put

(13)
$$H_1 = [C_1 \otimes D_{01}, \sigma_1(C_1) \otimes D_{11}, \sigma_2(C_1) \otimes D_{21}] \quad \text{and}$$
$$H_2 = [C_2 \otimes D_{02}, \sigma_1(C_2) \otimes D_{12}, \sigma_2(C_2) \otimes D_{22}].$$

THEOREM 9. *For $H_1$ and $H_2$ in (13), we have:*

(i) *the matrix $H_1$ is an $\mathrm{OA}(n_1 s_{11} s_{21}, (s_{11} s_{21})^{k_0} s_{11}^{k_1} s_{21}^{k_2})$;*

(ii) *the matrix $H_2$ becomes an $\mathrm{OA}(n_2 s_{12} s_{22}, (s_{12} s_{22})^{k_0} s_{12}^{k_1} s_{22}^{k_2})$ after the levels of the $s_{11} s_{21}$-level factors are collapsed according to $\delta_0$ and the levels of $s_{j1}$-level factors are collapsed according to $\delta_j$ for $j = 1, 2$.*

This theorem can be readily verified by following the theorem in Wang (1996) and the definition of $\delta_j$.

We now give a simple method for constructing a type of matrix $D = [D_0, D_1, D_2]$ required by the preceding theorem. (As a side note, this construction is related to Problem 6.17 in HSS, page 144.) For $j = 1, 2$, take $D_j$ to be a $D(b_j, c_j, s_j)$ based on $\mathrm{GF}(s_j)$. Let $c_0$ be any integer between 1 and $\min(c_1, c_2)$. For $j = 1, 2$, partition $D_j$ as $D_j = [D_{j0}, D_{j1}]$, where both $D_{10}$ and $D_{20}$ have $c_0$ columns. Let $\alpha_{i,k}$ denote the $(i, k)$th entry of $D_1$ and $\beta_{j,k}$ the $(j, k)$th entry of $D_2$. Construct a $b_1 b_2 \times c_0$ matrix $D_0$ whose entry in the $((i-1)b_2 + j)$th row and $k$th column is $(\alpha_{i,k}, \beta_{j,k})$, $1 \le i \le b_1, 1 \le j \le b_2, 1 \le k \le c_0$. Define $D = [D_0, D_{11}^*, D_{21}^*]$, where the $((i-1)b_2 + j)$th row of $D_{11}^*$ is the $i$th row of $D_{11}$ and the $((i-1)b_2 + j)$th row of $D_{21}^*$ is the $j$th row of $D_{21}$.

LEMMA 7. *For $D$ constructed above, we have:*

(i) *the matrix $D_0$ is a difference matrix $D(b_1 b_2, c_0, s_1 s_2)$;*
(ii) *for $j = 1, 2$, the matrix $D_{j1}^*$ is a difference matrix $D(b_1 b_2, c_j - c_0, s_j)$;*
(iii) *for $j = 1, 2$, $[\sigma_j(D_0), D_{j1}^*]$ is a difference matrix $D(b_1 b_2, c_j, s_j)$.*

EXAMPLE 13 [An $\mathrm{OA}(144, 12^2 4^2 3^1)$ containing an $\mathrm{OA}(72, 6^2 2^2 3^1)$]. Let $C_{10}$ be the vector listing all level combinations of $\mathrm{GF}(4)$ and $\mathrm{GF}(3)$. Denote by $0, 1, x, x+1$ the elements of $\mathrm{GF}(4)$ are $0, 1, 2$ the elements of $\mathrm{GF}(3)$. For $j = 1, 2$, let $C_{1j}$ be the column vector listing all $j$th digits of $C_{10}$. The transpose of the matrix $C_1 = (C_{10}, C_{11}, C_{12})$ is given by

$$\begin{bmatrix} 00 & 01 & 02 & 10 & 11 & 12 & x0 & x1 & x2 & (x+1)0 & (x+1)1 & (x+1)2 \\ 0 & 0 & 0 & 1 & 1 & 1 & x & x & x & x+1 & x+1 & x+1 \\ 0 & 1 & 2 & 0 & 1 & 2 & 0 & 1 & 2 & 0 & 1 & 2 \end{bmatrix}.$$

Let $C_2 = (C_{20}, C_{21}, C_{22})$ be the submatrix of $C_1$ consisting of the first six rows. Take $D_1$ to be the following $D(4,4,4)$

$$\begin{bmatrix} 0 & 0 & 0 & 0 \\ 0 & 1 & x & x+1 \\ 0 & x & x+1 & 1 \\ 0 & x+1 & 1 & x \end{bmatrix}$$

and $D_2$ to be the following $D(3,3,3)$

$$\begin{bmatrix} 0 & 0 & 0 \\ 0 & 1 & 2 \\ 0 & 2 & 1 \end{bmatrix}.$$

From Lemma 7, $D = [D_0, D_1, D_2]$ is

$$\begin{bmatrix} 00 & 00 & 0 & 0 & 0 \\ 00 & 01 & 0 & 0 & 2 \\ 00 & 02 & 0 & 0 & 1 \\ 00 & 10 & x & x+1 & 0 \\ 00 & 11 & x & x+1 & 2 \\ 00 & 12 & x & x+1 & 1 \\ 00 & x0 & x+1 & 1 & 0 \\ 00 & x1 & x+1 & 1 & 2 \\ 00 & x2 & x+1 & 1 & 1 \\ 00 & (x+1)0 & 1 & x & 0 \\ 00 & (x+1)1 & 1 & x & 2 \\ 00 & (x+1)2 & 1 & x & 1 \end{bmatrix}.$$

Put $H_1 = [C_{10} \otimes D_0, C_{11} \otimes D_1, C_{12} \otimes D_2]$ and $H_2 = [C_{20} \otimes D_0, C_{21} \otimes D_1, C_{22} \otimes D_2]$. From Theorem 9, $H_1$ is an OA$(144, 12^2 4^2 3^1)$ and $H_2$ becomes an OA$(72, 6^2 2^2 3^1)$ after the following level collapsing: for the 12-level factors, using $\{00, x0\} \to 00$, $\{10, (x+1)0\} \to 10$, $\{01, x1\} \to 01$, $\{11, (x+1)1\} \to 11$, $\{02, x2\} \to 02$, $\{12, (x+1)2\} \to 12$; for the 4-level factors, using $\{0, 2\} \to 0$ and $\{1, 3\} \to 1$.

**8. Generation of nested space-filling designs.** In this section, we discuss the problem of using NOAs to generate NSFDs. Throughout, we assume the factors are quantitative and each of them takes values in the interval $[0, 1]$. When we say that a design is space-filling or achieves uniformity in low dimensions, we mean that, when projected onto low dimensions, the design points are evenly scattered in the design region. For this problem, we present an approach following the procedure in QTW used for the same problem. Unlike QTW, the present approach covers both NOAs with equal levels and with mixed levels. Consider an NOA$(H_1, H_2)$, where $H_1$ is an OA$(n_1, s_{11}^{\gamma_1} \cdots s_{k1}^{\gamma_k})$ with $m = \sum_{i=1}^{k} \gamma_i$, $H_2$ is a submatrix of $H_1$ and $H_2$ becomes an OA$(n_2, s_{12}^{\gamma_1} \cdots s_{k2}^{\gamma_k})$ after the levels of the $s_{j1}$-level factors are collapsed into $s_{j2}$ levels according to a projection $\delta_j$. If $k = 1$, this array reduces to an NOA with equal levels.

The first step in constructing an OA-based Latin hypercube design $D_l$ using $H_1$ is to relabel the $s_{j1}$ levels of $H_1$, currently represented by the elements of a Galois field (or other mathematical structures), as $1, \ldots, s_{j1}$. Note that the projection $\delta_j$ divides the $s_{j1}$ levels into $s_{j2}$ groups, each of size $e_j = s_{j1}/s_{j2}$, and two levels belong to the same group if their projected values match. To ensure that the subset of $D_l$ corresponding to $H_2$ has good space-filling properties, we label the $s_{j1}$ levels of any $s_{j1}$-level factor in $H_1$ in such a way that the group of levels that are mapped to the same level should form a consecutive subset of $\{1, \ldots, s_{j1}\}$. The $s_{j2}$ groups are arbitrarily labeled as groups $1, \ldots, s_{j2}$, and the $e_j$ levels within the $i$th group are arbitrarily labeled as $(i-1)e_j + 1, \ldots, (i-1)e_j + e_j$ for $i = 1, \ldots, s_{j2}$.

After labeling the levels of the $s_j$-level factors of $H_1$ as $1, \ldots, s_{j1}$, $j = 1, \ldots, k$, as discussed above, we now use this array to obtain an OA-based Latin hypercube design as described in Section 2.1. Let $D_l$ denote the set of points and $D_h$ be the subset of $D_l$ corresponding to $H_2$. Then (i) $D_l$ achieves maximum uniformity in one dimension and, when $D_l$ is projected onto the dimensions of an $s_{j1}$-level factor and an $s_{k1}$-level factor, the points achieve uniformity on $s_{j1} \times s_{k1}$ grids; and (ii) $D_h$ is a subset of $D_l$ and, when $D_l$ is projected onto the dimensions of an $s_{j1}$-level factor and an $s_{k1}$-level factor, the points achieve uniformity on $s_{j2} \times s_{k2}$ grids.

An example is given to illustrate the above procedure.

EXAMPLE 14. Consider the NOA in Example 8, where $H_1$ is an OA(64, 4, 8), $H_2$ is a submatrix of $H_1$ and $\phi(H_2)$ is an OA(32, 4, 4). The four groups of levels of $H_1$ are $\{0, x^2\}$, $\{1, x^2 + 1\}$, $\{x, x^2 + x\}$ and $\{x + 1, x^2 + x + 1\}$. We label $\{0, x^2\}$ as levels 1 and 2, $\{1, x^2 + 1\}$ as levels 3 and 4, $\{x, x^2 + x\}$ as levels 5 and 6 and $\{x + 1, x^2 + x + 1\}$ as levels 7 and 8. Table 4 presents the array $H_1$ after using such labeling, where $H_2$ correspond to runs 1, 2, 7–10, 15–18, 23–26, 31–34, 39–42, 47–50, 55–58, 63–64. We then use $H_1$ to construct an OA-based Latin hypercube design $D_l$ for $x_1$ to $x_4$. Now choose $D_h$ to be the subset of $D_l$ corresponding to $H_2$. The points in any bivariate projection of $D_h$ achieve uniformity on the $4 \times 4$ grids in two dimensions. The points in the bivariate projections of $D_l$ also achieve similar uniformity.

**9. Discussions and concluding remarks.** Multiple computer experiments with different levels of accuracy have become prevalent in business, engineering and science for studying complex real world systems. NSFDs are attractive for such experiments. Several methods are proposed for constructing various families of NOAs, which can be used to generate many new NSFDs. In the development of these methods, two new discrete mathematics concepts, called nested orthogonal arrays and nested difference matrices, are introduced. These concepts should be further studied in their own right.

NSFDs can also be used in validation of computer models, that is, testing the accuracy of a computer model against some field data [Bayarri et al. (2007), Kennedy and O'Hagan (2001) and Oberkampf and Trucano (2007)]. Let $D_c$ denote the set of design points for the computer model and $D_f$ denote the set of design points for the corresponding physical experiment used as a benchmark in the validation. Unlike the situation of $D_l$ and $D_h$, $D_c$ should have more columns than $D_f$ because of the need of accommodating calibration (tuning) parameters that appear in the computer model only. Precisely, construction of $D_c$ and $D_f$ is guided by the following requirements:

TABLE 4
*The $H_1$ matrix in Example [14](#)*

| Run # | $x_1$ | $x_2$ | $x_3$ | $x_4$ | Run # | $x_1$ | $x_2$ | $x_3$ | $x_4$ |
|---|---|---|---|---|---|---|---|---|---|
| 1 | 1 | 1 | 1 | 1 | 33 | 2 | 2 | 2 | 2 |
| 2 | 1 | 3 | 5 | 7 | 34 | 2 | 4 | 6 | 8 |
| 3 | 1 | 2 | 7 | 8 | 35 | 2 | 1 | 8 | 7 |
| 4 | 1 | 4 | 3 | 2 | 36 | 2 | 3 | 4 | 1 |
| 5 | 1 | 5 | 2 | 6 | 37 | 2 | 6 | 1 | 5 |
| 6 | 1 | 7 | 6 | 4 | 38 | 2 | 8 | 5 | 3 |
| 7 | 1 | 6 | 8 | 3 | 39 | 2 | 5 | 7 | 4 |
| 8 | 1 | 8 | 4 | 5 | 40 | 2 | 7 | 3 | 6 |
| 9 | 3 | 3 | 3 | 3 | 41 | 4 | 4 | 4 | 4 |
| 10 | 3 | 1 | 7 | 5 | 42 | 4 | 2 | 8 | 6 |
| 11 | 3 | 4 | 5 | 6 | 43 | 4 | 3 | 6 | 5 |
| 12 | 3 | 2 | 1 | 4 | 44 | 4 | 1 | 2 | 3 |
| 13 | 3 | 7 | 4 | 8 | 45 | 4 | 8 | 3 | 7 |
| 14 | 3 | 5 | 8 | 2 | 46 | 4 | 6 | 7 | 1 |
| 15 | 3 | 8 | 6 | 1 | 47 | 4 | 7 | 5 | 2 |
| 16 | 3 | 6 | 2 | 7 | 48 | 4 | 5 | 1 | 8 |
| 17 | 5 | 5 | 5 | 5 | 49 | 6 | 6 | 6 | 6 |
| 18 | 5 | 7 | 1 | 3 | 50 | 6 | 8 | 2 | 4 |
| 19 | 5 | 6 | 3 | 4 | 51 | 6 | 5 | 4 | 3 |
| 20 | 5 | 8 | 7 | 6 | 52 | 6 | 7 | 8 | 5 |
| 21 | 5 | 1 | 6 | 2 | 53 | 6 | 2 | 5 | 1 |
| 22 | 5 | 3 | 2 | 8 | 54 | 6 | 4 | 1 | 7 |
| 23 | 5 | 2 | 4 | 7 | 55 | 6 | 1 | 3 | 8 |
| 24 | 5 | 4 | 8 | 1 | 56 | 6 | 3 | 7 | 2 |
| 25 | 7 | 7 | 7 | 7 | 57 | 8 | 8 | 8 | 8 |
| 26 | 7 | 5 | 3 | 1 | 58 | 8 | 6 | 4 | 2 |
| 27 | 7 | 8 | 1 | 2 | 59 | 8 | 7 | 2 | 1 |
| 28 | 7 | 6 | 5 | 8 | 60 | 8 | 5 | 6 | 7 |
| 29 | 7 | 3 | 8 | 4 | 61 | 8 | 4 | 7 | 3 |
| 30 | 7 | 1 | 4 | 6 | 62 | 8 | 2 | 3 | 5 |
| 31 | 7 | 4 | 2 | 5 | 63 | 8 | 3 | 1 | 6 |
| 32 | 7 | 2 | 6 | 3 | 64 | 8 | 1 | 5 | 4 |

(i) $D_c$ contains all factors of $D_f$ and has additional columns to accommodate the calibration parameters.

(ii) When restricted to the shared factors, $D_f \subset D_c$.

(iii) Both $D_c$ and $D_f$ have good space-filling properties.

With slight modifications, our construction methods for $D_h$ and $D_l$ can give $D_f$ and $D_c$ that satisfy the above requirements. For illustration, we modify the construction in Section 3.1, where $F = \mathrm{GF}(2^{m+1})$, $G = \mathrm{GF}(2^m)$, $m \geq 2$, $D_0$ is a $D(2^{m+1}, 2^{m+1}, 2^{m+1})$ and $D_1$ is a $D(2^{m+1}, 2^4, 2^{m+1})$. Take $D_1^*$ to be $D_0$. Then $D_1^*$ has more columns than $D_1$. Next replace $D_1$ by $D_1^*$ and follow through the steps in Section 3.1 and the construction in (10). Let $A$ be the $\mathrm{OA}(n, c, s_1)$ used in Theorem 4. Then we have the following results: for $m \geq 2$,

(i) the matrix $D_1^*$ is a $D(2^{m+1}, 2^{m+1}, 2^{m+1})$;
(ii) the matrix $\phi(D_2)$ is a $D(2^m, 2^2, 2^m)$;
(iii) the matrix $H_1^* = A \otimes D_1^*$ is an $\mathrm{OA}(n2^{m+1}, 2^{m+1}c, 2^{m+1})$;
(iv) for the shared $4m_2$ factors, $H_2 = A \otimes D_2$ is a submatrix of $H_1^*$ and $\delta(H_2)$ is an $\mathrm{OA}(n2^m, 4c, 2^m)$.

As in Section 8, we use $(H_1^*, H_2)$ to generate a pair of nested designs for $D_c$ and $D_f$, where $D_c$ has $2^{m+1}c$ columns and $D_f$ has $4c$ columns and both have good space-filling properties.

Extensions of the present work can be made in several directions. First, similar to the construction of OA-based Latin hypercubes designs [Tang (1993, 1994) and Leary, Bhaskar and Keane (2003)], it is possible to produce multiple NSFDs based on a given NOA. In a separate article, we plan to use both distance and correlation criteria to construct optimal NSFDs. Second, the constructed NOAs in this article have strength 2 that can guarantee uniformity in two dimensions only. The proposed methods can be extended to produce NOAs with higher-dimensional stratification by exploring nesting in difference matrices with strength 3 or higher [Hedayat, Stufken and Su (1996)]. Another possibility is to use quasi-Monte Carlo sequences, like nets [Niederreiter (1992)]. A paper in preparation will address the issue of constructing nested nets. Third, it is worth studying the sampling properties of NSFDs. Fourth, a natural extension of the present work is to construct NSFDs for experiments with more than two levels of accuracy. One way to achieve this is to extend the method in QTW to directly construct NOAs with more sophisticated nesting, that is, a 32-run OA contains a 16-run OA that contains an 8-run OA. Another possibility is to modify the method in Section 3 to obtain NDMs with nesting at more than two levels and then use them to produce the desired NOAs. Finally, given the close connections between OAs and coding theory, it should be possible to use coding-theoretical techniques to construct new NOAs. We are currently exploring this issue.

## APPENDIX

$D(12, 12, 4)$ from Seberry (1979)

| 00 | 00 | 00 | 00 | 00 | 00 | 00 | 00 | 00 | 00 | 00 | 00 |
|----|----|----|----|----|----|----|----|----|----|----|----|
| 00 | 00 | 00 | 01 | 01 | 01 | 11 | 11 | 11 | 10 | 10 | 10 |
| 00 | 00 | 00 | 11 | 11 | 11 | 10 | 10 | 10 | 01 | 01 | 01 |
| 00 | 11 | 01 | 10 | 01 | 11 | 01 | 10 | 00 | 11 | 00 | 00 |
| 00 | 11 | 01 | 11 | 10 | 01 | 00 | 01 | 10 | 10 | 11 | 00 |
| 00 | 11 | 01 | 01 | 11 | 10 | 10 | 00 | 01 | 00 | 10 | 11 |
| 00 | 01 | 10 | 11 | 00 | 10 | 01 | 00 | 11 | 01 | 11 | 10 |
| 00 | 01 | 10 | 10 | 11 | 00 | 11 | 01 | 00 | 10 | 01 | 11 |
| 00 | 01 | 10 | 00 | 10 | 11 | 00 | 11 | 01 | 11 | 10 | 01 |
| 00 | 10 | 11 | 01 | 10 | 00 | 01 | 11 | 10 | 01 | 00 | 11 |
| 00 | 10 | 11 | 00 | 01 | 10 | 10 | 01 | 11 | 11 | 01 | 00 |
| 00 | 10 | 11 | 10 | 00 | 01 | 11 | 10 | 01 | 00 | 11 | 01 |

$D(12, 6, 12)$ based on $(\mathbb{Z}_2 \oplus \mathbb{Z}_6, +)$ [Dulmage, Johnson and Mendelsohn (1961)]

| 00 | 00 | 00 | 00 | 00 | 00 |
|----|----|----|----|----|----|
| 00 | 01 | 03 | 12 | 04 | 10 |
| 00 | 02 | 10 | 01 | 15 | 12 |
| 00 | 03 | 01 | 15 | 14 | 02 |
| 00 | 04 | 13 | 05 | 02 | 11 |
| 00 | 05 | 15 | 13 | 11 | 01 |
| 00 | 10 | 02 | 03 | 12 | 13 |
| 00 | 11 | 12 | 14 | 10 | 15 |
| 00 | 12 | 05 | 02 | 13 | 04 |
| 00 | 13 | 04 | 11 | 01 | 14 |
| 00 | 14 | 11 | 10 | 03 | 05 |
| 00 | 15 | 14 | 04 | 05 | 03 |

**Acknowledgments.** The authors thank the Editor, the Associate Editors, and two referees for their valuable comments and suggestions that improved the presentation of this article.

P. Z. G. Qian  
Department of Statistics  
University of Wisconsin–Madison  
Madison, Wisconsin 53706  
USA  
E-mail: peterq@stat.wisc.edu

M. Ai  
LMAM, School of Mathematical Sciences  
Peking University  
Beijing 100871  
China  
E-mail: myai@math.pku.edu.cn

C. F. J. Wu  
H. Milton Stewart School of Industrial  
 and Systems Engineering  
Georgia Institute of Technology  
755 Ferst Drive NW, Atlanta  
Georgia 30332-0205  
USA  
E-mail: jeffwu@isye.gatech.edu